\documentclass[12pt,english]{article}
\usepackage{theorem}
\usepackage{a4wide}
\usepackage{graphics}
\usepackage{epsfig}
\usepackage{amssymb}
\usepackage{latexsym}
\usepackage{amsmath}

\newtheorem{thh}{Theorem}[section]
\newtheorem{df}[thh]{Definition}

\newtheorem{lem}[thh]{Lemma}
\newtheorem{cor}[thh]{Corollary}
\newtheorem{prop}[thh]{Proposition}

\title{Families of $k$-derivations on $k$-algebras}
\author{Philippe Bonnet}
\date{}

\newcommand{\dem}{{\em Proof: }}
\newcommand{\qed}{\begin{flushright} $\blacksquare$\end{flushright}}
\newcommand{\der}{ \partial}
\newcommand{\MM}{{\cal{M}}}

\newcommand{\FF}{{\cal{F}}}
\newcommand{\CC}{\mathbb C}

\newcommand{\PR}{{\cal{P}}}

\begin{document}
\maketitle

\begin{center} { \small
Mathematisches Institut, Universit\"at Basel\\
Rheinsprung 21, 4051 Basel, Switzerland\\
bonnet@math-lab.unibas.ch}
\end{center}

\begin{abstract}
Let $A$ be an integral $k$-algebra of finite type over a field $k$ of characteristic
zero. Let $\FF$ be a family of $k$-derivations on $A$
and $M_{\FF}$ the $A$-module spanned by $\FF$. In this paper, we generalise a result due to A. Nowicki
and construct
an element $\der$ of $M_{\FF}$ such that $ker \; \der =
\cap_{d \in \FF} ker \; d$. Such a derivation is called $\FF$-minimal.
Then we establish a density theorem for $\FF$-minimal derivations in $M_{\FF}$.
\end{abstract}

\section{Introduction}

Let $A$ be an integral $k$-algebra of finite type over a field $k$ of characteristic
zero.
For convenience, we will say that a sub-algebra $B$ of $A$
is {\em algebraically closed in $A$} if every element $a$ of $A$ that is algebraic
over $B$ belongs to $B$. Let $\FF$ be a family of $k$-derivations on $A$.
In this paper, we are interested in describing the kernel of this family, i.e.
the following set:
$$
ker \; \FF = \cap_{d \in \FF} \; ker \; d
$$
Let $M_{\FF}$ be the $A$-module spanned by the elements of $\FF$. By analogy
with the theory of foliations, we say that an element $f$ of $A$ is a
{\em first integral of a $k$-derivation
$d$} if $d(f)=0$ and $f\not\in k$.
Similarly $f$ is a {\em first integral of $\FF$} if
$d(f)=0$ for every $d\in \FF$, and $f\not\in k$. First integrals correspond
to the notion of constants for a derivation (see \cite{No}), except that
they must not belong to the coefficient field $k$.

The description of the kernels $ker \; \FF$ is usually quite tricky because of their complexity.
Indeed, since Nagata's works (see \cite{Na}), it is well-known that the sub-algebra
$B=ker\; \FF$ neednot be finitely generated. Nagata's construction uses locally nilpotent
derivations on a $k$-algebra $A$ of Krull dimension $n\geq 32$. This result has been refined
by Deveney and Finston, who constructed a locally nilpotent $k$-derivation on $k[x_1,...,x_7]$
whose kernel is not finitely generated (see \cite{De-F}). Recently this result has been improved by
Daigle and Freudenburg (see \cite{Da-F}), with an example of a locally nilpotent derivation
on $k[x_1,...,x_5]$ having as kernel a non-finitely generated algebra. In contrast, such
behaviours do not occur in low dimensions. For instance, derivations on $k[x_1,...,x_n]$ have as
kernel a finitely generated $k$-algebra if $n\leq 3$ (see \cite{Na2}).

In what follows, we will choose to express $ker\; \FF$ not as a $k$-algebra, but in terms
of the derivations involved in its construction.
Our starting point is an article of Nowicki (see \cite{No}) where he proved the
following two theorems.

\begin{thh} {\rm{(\cite{No})}} \label{Nowicki1}
Let $A$ be an integral $k$-algebra of finite type over a field $k$ of characteristic
zero. Let $\FF$ be a family of $k$-derivations on $A$.
Then there exists a $k$-derivation $d$ on $A$ such that $ker \; d=ker\; \FF$.
\end{thh}

\begin{thh} {\rm{(\cite{No})}} \label{Nowicki2}
Let $A$ be an integral $k$-algebra of finite type over a field $k$ of characteristic
zero. Let $B$ be a sub-algebra of $A$. Then $B$ is
algebraically closed in $A$ if and only if $B$ is the kernel of a $k$-derivation
on $A$.
\end{thh}
The proof of theorem \ref{Nowicki1} is very elegant and uses Noether normalisation
lemma. However the construction of the derivation $d$ is independent of the family
$\FF$, and it only uses the fact that the ring $B=ker \; \FF$ is algebraically
closed in $A$. In this paper we will refine theorem \ref{Nowicki1}, and we
will express the derivation $d$ in terms of the elements of $\FF$.
More precisely:

\begin{thh} \label{principal}
Let $A$ be an integral $k$-algebra of finite type over a field $k$ of characteristic
zero. Let $\FF=\{d_i\}_{i\in I}$ be a family of
$k$-derivations on $A$. Then there exists an $A$-linear combination
$d=\sum a_id_i$ such that $
ker \; d = ker \; \FF$.
\end{thh}
A $k$-derivation $d$ in $M_{\FF}$ is {\em $\FF$-minimal} if
$ker \; d= ker \; \FF$. Here $\FF$-minimality means that
the kernel of $d$ is smallest among all kernels of $k$-derivations
in $M_{\FF}$.
We denote by $M_{\FF, min}$ the set of all $\FF$-minimal $k$-derivations on
$A$. In terms of first integrals, Theorem \ref{principal} can be reinterpreted
as follows:

\begin{cor} \label{cons1}
Let $A$ be an integral $k$-algebra of finite type over a field $k$ of characteristic
zero. Let $\FF=\{d_i\}_{i\in I}$ be a family of
$k$-derivations on $A$. If every $A$-linear combination $d =\sum a_id_i$
admits a first integral, then $\FF$ admits a first integral.
\end{cor}
Let $d =\sum a_id_i$ be an $A$-linear combination of elements of $\FF$.
A priori its first integrals (if any) depend on the coefficients $a_i$. But if
every such combination has a first integral, then the previous corollary
asserts that we can choose a first integral $f$ that is independent of the
$a_i$. In particular $d_i(f)=0$ for any $i\in I$.

\begin{df}
Let $k$ be an algebraically closed field of characteristic zero. Let
$E$ be a $k$-vector space and $\Omega$ a subset of $E$. The set $\Omega$
is residual in $E$ if, for any finite dimensional $k$-subspace $F$ of $E$,
$\Omega \cap F$ is a countable intersection of Zariski open sets
of $F$ (possibly empty).
\end{df}
Note that if $k=\CC$ and $\Omega \cap F\not=\emptyset$, then
$\Omega \cap F$ is dense in $F$ for the Zariski and metric
topologies on $F$. This latter assertion is based on Baire's Theorem
about countable intersection of dense open sets in a complete space,
and also on the fact that every non-empty Zariski open set is dense
in $F$ for the metric topology. So residuality can be seen as a version of
density adapted to infinite dimensional spaces. With this definition,
theorem \ref{principal} yields the following
result:

\begin{thh} \label{cons2}
Let $k$ be an algebraically closed field of characteristic zero. Let $A$ be an
integral
$k$-algebra of finite type. Let $\FF$ be a family of
$k$-derivations on $A$. Then $M_{\FF, min}$ is a non-empty residual subset
of $M_{\FF}$.
\end{thh}
At the end of this paper, we will give a proof of theorem \ref{Nowicki2}
based on theorem \ref{principal}, and we will illustrate the notion of residuality
with an example. Note that in this paper, we have not investigated the field
of rational first integrals, i.e. the elements $f$ of the fraction field $K(A)$
of $A$ such that $d(f)=0$ for any $d$ in $\FF$. This field can be extremely large
compared to the kernel of $\FF$, as is the case for the Euler vector field:
$$
d=x_1\frac{\der}{\der x_1} + ...+ x_n\frac{\der}{\der x_n}
$$
on $k^n$, but it can also be reduced to the field $k$ (see for instance \cite{No2}
or \cite{Jou}). One question could be to try and find an analogue to theorem \ref{principal}.
More precisely, given a collection of $k$-derivations $\{d_i\}$ on a field $K$ of finite
transcendence degree, with $L=\cap ker d_i|K$, does there exist a $K$-linear combination
$d$ of the $d_i$ having $L$ as its kernel? 

\section{Reduction to a couple of derivations}

Throughout this paper, $A$ will be an integral $k$-algebra of finite type over a
field $k$ of
characteristic zero. Let $\FF$ be a family of
$k$-derivations on $A$. We will say that {\em the triplet $(k,A,\FF)$ enjoys
the property ${\cal{P}}$} if there exists an A-linear combination
$d=\sum a_i d_i$ of elements $d_i$ of $\FF$
such that:
$$
ker \; \der = ker \; \FF
$$
In this section, we are going to see how to restrict the proof of theorem \ref{principal}
to the case of a couple of derivations enjoying some remarkable properties. More precisely:

\begin{prop} \label{reduction}
$\PR$ holds for any triplet $(k,A,\FF)$ if and only if $\PR$ holds for any triplet
$(k',A',\FF')$, where $\FF'=\{d' _1,d' _2\}$ is a couple of $k'$-derivations
satisfying the two conditions: $(1)$ $ker\; d' _1 \cap ker\;
d' _2=k'$ and $(2)$ there exist two elements $x_1,x_2$ of $A'$
such that $d' _i(x_i)=1$ and $d' _i(x_j)=0$
if $i\not=j$.
\end{prop}
The proof of this proposition is a consequence of the following lemmas.
 
\begin{lem}
$\PR$ holds for any triplet $(k,A,\FF)$ if and only if it holds for any $(k',A',\FF')$,
where $\FF'$ is a finite family.
\end{lem}
\dem One direction is clear. For the other, consider the triplet $(k,A,\FF)$ and
let $Der_k(A,A)$ be the $A$-module of $k$-derivations on $A$. Since $A$ is a finite
$k$-algebra, $Der_k(A,A)$ is a noetherian $A$-module. So the $A$-submodule $M_{\FF}$ of
$Der_k(A,A)$ spanned by the elements of $\FF$ is finitely generated. Let $\FF'=\{d_1,..,d_r\}$ be a
finite subset of $\FF$ whose elements span the $A$-module $M_{\FF}$. Then we have:
$$
ker \FF = ker \; \FF' =ker\; d_1 \cap ...\cap ker \; d_r
$$
The inclusion $ker \FF \subseteq ker \; \FF'$ is obvious. Conversely let $f$ be an element of $A$
such that $d_1(f)=...=d_r(f)=0$. For any $k$-derivation $d$ of $\FF$, there exist some
elements $a_1,...,a_r$ of $A$ such that:
$$
d= a_1d_1 +...+a_rd_r
$$
Therefore $d(f)=a_1d_1(f)+...+a_rd_r(f)=0$, $f$ belongs to $ker\; \FF$ and $ker\; \FF=ker\;
\FF'$. Assume that $\PR$ holds for any $(k',A',\FF')$, where $\FF'$ is finite. Apply this
property to the triplet $(k,A,\FF')$. Then there exists an $A$-linear combination $\der$
of elements of $\FF'$ such that $ker \; \der= ker \; \FF'$. Since $ker \; \FF'=ker \; \FF$
and $\FF'$ is a subset of $\FF$, the result follows.
\qed

\begin{lem}
$\PR$ holds for any triplet $(k,A,\FF)$, where $\FF$ is finite, if and only if it holds
for any $(k',A',\FF')$,
where $\FF'=\{d_1,d_2\}$.
\end{lem}
\dem One direction is clear. The other will be proved by induction on the order $r$
of $\FF$. If $r=1$, then $\FF$ consists of one derivation $d_1$ and we choose $\der=d_1$.
If $r=2$, then $\FF$ is a couple of derivations and the result follows by assumption.
Assume the property holds to the order $r\geq 2$, and let $(k,A,\FF)$ be a triplet
such that $\FF=\{d_1,...d_{r+1}\}$. By assumption there exists a $k$-derivation
$\der'=a_1d_1+a_2d_2$, where every $a_i$ belongs to $A$, such that:
$$
ker \; \der' =ker\; d_1 \cap ker\; d_2
$$
Consider the family $\{\der',d_3,..,d_{r+1}\}$. Since the property holds
to the order $r$, there exists a $k$-derivation $\der=b'\der' + ...+b_{r+1}d_{r+1}$,
where $b',b_3,..,b_{r+1}$ belong to $A$, such that:
$$
ker \; \der = ker \; \der' \cap ker \; d_3 \cap ...\cap d_{r+1}= ker \; d_1 \cap ker \; d_2
\cap ...\cap ker \; d_{r+1}
$$
Since $\der'$ is a $A$-linear combination of $d_1,d_2$, $\der$ belongs to the $A$-module
spanned by $d_1,...,d_{r+1}$, and the result follows.
\qed 

\begin{lem}
$\PR$ holds for any triplet $(k,A,\FF)$, where $\FF=\{d_1,d_2\}$, if and only if it holds
for any $(k',A',\FF')$, where $\FF'=\{d' _1,d' _2\}$ is a couple of $k'$-derivations for
which there exist two elements $x_1,x_2$ of $A'$ such that $d' _i(x_i)=1$ and $d '_i(x_j)=0$
if $i\not=j$.
\end{lem}
\dem One direction is clear. For the other, assume that $\PR$ holds
for any $(k',A',\FF')$, where $\FF'=\{d' _1,d' _2\}$ is a couple of $k'$-derivations for
which there exist two elements $x_1,x_2$ of $A'$ such that $d' _i(x_i)=1$ and $d' _i(x_j)=0$
if $i\not=j$. Let $(k,A,\FF)$ be a triplet for which $\FF=\{d_1,d_2\}$. If $ker \; d_1
\subseteq ker \; d_2$ or $ker \; d_2 \subseteq ker \; d_1$, then we have:
$$
ker \; \FF = ker \; d_1 \quad \mbox{or} \quad ker \; \FF = ker \; d_2
$$
and $\PR$ holds by choosing either $\der =d_1$ or $\der=d_2$. So we may assume that:
$$
ker \; d_2 \not\subseteq ker \; d_1 \quad \mbox{and} \quad ker \; d_1 \not\subseteq ker \; d_2
$$
By assumption, there exist two elements $x_1,x_2$ of $A$ such that:
$$
d_1(x_1)\not=0, \quad d_2(x_1)=0, \quad d_2(x_2)\not=0, \quad d_1(x_2)=0
$$
We set $p=d_1(x_1)d_2(x_2)$ and consider the triplet:
$$(k',A',\FF')=\left(k, A[\frac{1}{p}], \left\{\frac{d_1}{d_1(x_1)},\frac{d_2}{d_2(x_2)}\right\} \right)
$$
By construction $A'$ is a $k'$-algebra of finite type, and it is a domain. Moreover $d'_1,d'
_2$ act on $A'$ as $k'$-derivations, $d' _i(x_i)=1$ and $d' _i(x_j)=0$ if $i\not=j$. So there
exists an $A'$-linear combination $\der'$ of $d'_1,d' _2$ such that:
$$
ker \; \FF' =ker \; \der'
$$
Up to replacing $\der'$ by $p^n \der'$ for $n$ big enough, we may assume that $\der'$ is an
$A$-linear combination of $d_1,d_2$. Let $\der$ be the restriction of $\der'$ to $A$. Let us
prove that $ker \; \FF =ker \; \der$.

The inclusion $ker \; \FF \subseteq ker \; \der$ is clear. Now let $x$ be an element of $A$
such that $\der(x)=0$. Then $\der'(x)=0$ in $A'$, and $x$ belongs to $ker\; d' _1 \cap
ker\; d' _2 \cap A$. Since $d_1$ (resp. $d_2$) is proportional to $d'_1$ (resp. $d' _2$), we
find:
$$
d_1(x)=d_2(x)=0
$$
Since $x$ belongs to $A$, $x$ belongs to $ker \; d_1 \cap ker\; d_2 =ker\;
\FF$ and the result follows.
\qed

\begin{lem}
The following assertions are equivalent:
\begin{itemize}
\item{$\PR$ holds for any triplet $(k,A,\FF)$, where $\FF=\{d_1,d_2\}$ is a couple of
$k$-derivations for which there exist two elements $x_1,x_2$ of $A$ such that
$d _i(x_i)=1$ and $d _i(x_j)=0$ if $i\not=j$.}
\item{$\PR$ holds for any triplet $(k',A',\FF')$, where $\FF'=\{d' _1,d' _2\}$ is a
couple of $k'$-derivations satisfying the two conditions: $(1)$ $ker\; d' _1 \cap ker\;
d' _2=k'$ and $(2)$ there exist two elements $x_1,x_2$ of $A'$ such that $d' _i(x_i)=1$ and $d' _i(x_j)=0$
if $i\not=j$.}
\end{itemize}
\end{lem}
\dem The first assertion implies clearly the second. Assume now that the second holds, and
let $(k,A,\FF)$ be a triplet satisfying the conditions of the first assertion. Let $k'$
be the fraction field of $ker\; \FF$, and consider the following triplet:
$$
(k',A',\FF')=\left(k', A\otimes_k k', \left \{d' _1,d' _2 \right\}\right)
$$
where every $d' _i$ acts on $A\otimes_k k'$ according to the following rule:
$$
d' _i(x\otimes f)=d_i (x)\otimes f
$$
If $S= ker \; \FF -\{0\}$, then $A'= A_S$ and $A'$ is a domain. Since $k'$ contains $k$,
$A'$ is a $k'$-algebra of finite type. Moreover $d' _i(a/b)=d_i(a)/b$ for any $a/b$
in $A_S$. In particular, every $d' _i$ is a $k'$-derivation. Let
$f=a/b$ be an element of $A_S$ such that $d' _1(f)=d' _2(f)=0$. By construction
we get $d_1(a)=d_2(a)=0$, $a$ belongs to $ker \; \FF$ and $f$ lies in $k'$.
Thus we have:
$$
ker \; \FF' =ker\; d' _1 \cap ker\; d' _2 =k'
$$
So the conditions of the second assertion hold, and there exists an $A_S$-linear
combination $\der'=a' _1 d' _1 + a' _2d' _2$ such that:
$$ker\; \der'=ker \; \FF' =k'$$
Up to a mutiplication by an element of $S$, we may assume that $a' _1$ and $a' _2$
belong to $A$. Denote by $\der$ the derivation $a' _1d_1 + a' _2d_2$. Let $f$
be an element of $ker \; \der$. Since $A$ is contained in $A_S$ and every $d' _i$
extends $d_i$ to $A_S$, we have $\der(f)=\der'(f)= 0$. So $d_i(f)=d' _i(f)=0$
for $i=1,2$ and $f$ belongs to $ker \; \FF$. Therefore we get:
$$
ker \; \der = ker \; \FF = ker\; d_1 \cap ker\; d_2 
$$
\qed

\section{Passage to a complete regular local ring} \label{redd}

Let $A$ be an integral $k$-algebra of finite type, and $\FF=\{d _1,d _2\}$ a couple
of $k$-derivations satisfying the conditions of proposition \ref{reduction}, i.e.
(1) $ker\; d _1 \cap ker\; d _2=k$ and (2) there exist two elements $x_1,x_2$ of $A$
such that $d _i(x_i)=1$ for all $i$ and $d_i(x_j)=0$ if $i\not=j$.
In
this section, we will see how to extend the $k$-derivations $d_1,d_2$
into a couple of $L$-derivations on a formal ring $L[[t_1,..,t_n]]$. This will
enable us to rewrite these derivations into a canonical form, that will
prove easier to handle.

\begin{prop} \label{complete}
Let $A$ be an integral $k$-algebra of finite type, and $\FF=\{d _1,d _2\}$ be a couple
of $k$-derivations satisfying the conditions of proposition \ref{reduction}. Let $x_1,x_2$ be
two elements of $A$ such that $d _i(x_i)=1$ for all $i$ and $d_i(x_j)=0$ if $i\not=j$. Then
there exist two elements $\lambda_1,\lambda_2$ of $k$, and an extension $L$ of $k$ such that:
\begin{itemize}
\item{$A$ is a subring of $L[[t_1,...,t_n]]$,}
\item{$d_1,d_2$ extend to $L$-derivations on $L[[t_1,...,t_n]]$,}
\item{$x_1 -\lambda_1,x_2-\lambda_2$ belong to the maximal ideal of $L[[t_1,...,t_n]]$.}
\end{itemize}
\end{prop}
The proof of this proposition will split into several lemmas.

\begin{lem}
Let $x_1,x_2$ be two elements of $A$ satisfying the conditions of proposition \ref{complete}.
Then $x_1,x_2$ are algebraically independent in the $k$-algebra $A$.
\end{lem}
\dem Assume there exists a non-zero polynomial $P$ in $k[u,v]$ such that
$P(x_1,x_2)=0$. We choose $P$ of minimal homogeneous degree with respect to $u,v$.
Since $d _i(x_i)=1$ for all $i$ and $d_i(x_j)=0$ if $i\not=j$, we get by derivation:
$$
d_1(P(x_1,x_2))=\frac{\partial P}{\partial u}(x_1,x_2)=0 \quad \mbox{and}
\quad d_2(P(x_1,x_2))=\frac{\partial P}{\partial v}(x_1,x_2)=0
$$
By minimality of the degree, this implies that $\frac{\partial P}{\partial u}=
\frac{\partial P}{\partial v}=0$. Therefore $P$ is constant and $P(x_1,x_2)=0$
implies that $P=0$, hence a contradiction.
\qed

\begin{lem}
There exist two elements $\lambda_1,\lambda_2$ of $k$, and a maximal ideal ${\cal{M}}$
of $A$ such that $x_1 -\lambda_1,x_2-\lambda_2$ belong to ${\cal{M}}$ and $A_{\cal{M}}$
is a regular local ring.
\end{lem}
\dem Up to localizing $A$ with respect to a non-zero element $g$ of $A$, we may assume
that $\Omega_{A_{(g)}/k}\simeq (\Omega_{A/k})_{(g)}$ is free. By generic smoothness (see \cite{Ei}),
$A'=A_{(g)}$ is regular over every maximal ideal ${\cal{M}}'$ not containing $g$.
Since $x_1,x_2$ are algebraically independent, the inclusion induces
an injective
map:
$$
L: k[u,v] \longrightarrow A', \quad P \longmapsto P(x_1,x_2)
$$
Therefore the map $L^*: Spec(A')\rightarrow Spec(k[u,v])$ is dominant, and there
exists an element $f\not=0$ of $k[u,v]$ such that every fibre $L^{-1}({\cal{P}})$
is non-empty for any maximal ideal ${\cal{P}}$ of $k[u,v]$ not containing
$f$. Since $f$ is non-zero, there exists a couple $(\lambda_1,\lambda_2)$
in $k^2$ such that $f(\lambda_1,\lambda_2)\not=0$. Consider the ideal:
$$
{\cal{P}}=(u -\lambda_1,v -\lambda_2)
$$
By construction, ${\cal{P}}$ is maximal in $k[u,v]$ and does not contain
$f$. So the fibre $L^{-1}({\cal{P}})$ is not empty. In particular, it contains
a maximal ideal ${\cal{M}}'$ of $A'$. If ${\cal{M}}$ denotes the intersection
${\cal{M}}'\cap A$, then ${\cal{M}}$ is a maximal ideal not containing $g$,
and we have the isomorphism
of $k$-algebras:
$$
A'_{{\cal{M}}'} \simeq A_{{\cal{M}}}
$$
Since $A'$ is regular over every maximal ideal, $A_{{\cal{M}}}$ is a regular local
ring. By construction ${\cal{M}}$ contains $x_1 -\lambda_1,x_2 -\lambda_2$, and the result
follows.
\qed
Since $\lambda_1,\lambda_2$ are annihilated by the $d_i$, we may replace $x_i$
by $x_i -\lambda_i$ without changing the conditions at the beginning of this section.
So we may assume that {\em $x_1,x_2$ belong to a maximal ideal ${\cal{M}}$ of $A$
such that
the $k$-algebra $A_{\cal{M}}$ is regular}. By an easy computation, we get
for any positive integer $r$:
$$
d_1({\cal{M}}^r)\subseteq {\cal{M}}^{r-1} \quad \mbox{and} \quad
d_2({\cal{M}}^r)\subseteq {\cal{M}}^{r-1}
$$
Therefore $d_1$ and $d_2$ are continuous on $A_{\cal{M}}$ for the ${\cal{M}}$-adic
topology, and they uniquely extend into a couple of $k$-derivations on
the ${\cal{M}}$-adic completion $R$ of $A_{\cal{M}}$. We still denote by $d_i$
this extension. Since $A_{\cal{M}}$
is regular and contains the field $k$, by Cohen Structure Theorem (see \cite{Ei}),
there exists an extension $L$ of $k$ such
that:
$$
R\simeq L[[t_1,..,t_n]]
$$
where $n$ is the Krull dimension of $A$. So the $d_i$ can be viewed as
$k$-derivations on $L[[t_1,..,t_n]]$. In order to get proposition \ref{complete},
we only need to check that: 

\begin{lem}
$d_1$ and $d_2$ are $L$-derivations on
$L[[t_1,..,t_n]]$. 
\end{lem}
\dem It suffices to prove that $d_i(L)=\{0\}$ for $i=1,2$. First note that
$L$ is isomorphic to $A/\cal{M}$. Since $A$ is a
finite $k$-algebra, the field $L$ is also a finite $k$-algebra.
So $L$ is a finite extension of $k$ (see \cite{Hum}). Let $\zeta$
be an element of $L$, and let $P$ be a polynomial in $k[t]$ of minimal
degree such that $P(\zeta)=0$. By derivation, we get:
$$
d_i(P(\zeta))=P'(\zeta)d_i(\zeta)=0
$$
By minimality of the degree, $P'(\zeta)\not=0$ and $d_i(\zeta)=0$. Since this
holds for any $\zeta$ in $L$, the result follows.
\qed

\section{Canonical form for a couple of $L$-derivations} 

In this section, we consider a couple of $L$-derivations $d_1,d_2$ on
$L[[t_1,..,t_n]]$, satisfying the following condition: there exist two
elements $x_1,x_2$ of $L[[t_1,..,t_n]]$ such that $d _i(x_i)=1$ for $i=1,2$
and $d_i(x_j)=0$ if $i\not=j$. We are going to search for a system of
parameters for which these derivations look simpler. Recall that a system
of parameters is a family of formal functions $s_1,...,s_n$ generating the
maximal ideal $(t_1,...,t_n)$ of $L[[t_1,..,t_n]]$.

\begin{lem} \label{champvecteur}
Let $d$ be an $L$-derivation on $L[[t_1,..,t_n]]$, where $n>1$. Assume there
exists a formal function $x_1$ such that $d(x_1)=1$. Then there exists
some formal functions $y_2,...,y_n$ such that $x_1,y_2,...,y_n$ is a system of
parameters and $d(y_i)=0$ for any $i$. In particular $d=\partial/\partial x_1$
in this system of parameters.
\end{lem}
\dem Let ${\cal{M}}$ be the maximal ideal of $L[[t_1,..,t_n]]$. Since $d$
is an $L$-derivation, we have $d(\MM^2)\subseteq \MM$. So $x_1$ belongs to
$\MM - \MM^2$ because $d(x_1)=1$. Let $y^0 _2,...,y^0 _n$ be a system of formal
functions such that $x_1,y^0 _2,...,y^0 _n$ form a basis of $\MM/\MM^2$.
By Nakayama Lemma, $x_1,y^0 _2,...,y^0 _n$ is a system of parameters. If
$\lambda_i=d(y_i)(0)$ and $y^1 _i= y^0 _i -\lambda_i x_1$, then we find:
$$
d(y^1 _i)(0)=\lambda_i -\lambda_id(x_1)(0)=0
$$
and $x_1,y^1 _2,...,y^1 _n$ is again a system of parameters. For any
integer $k>0$, we are going to construct a system $y^k _2,...,y^k _n$
of formal functions satisfying the following conditions:

\begin{itemize}
\item{$x_1,y^k _2,...,y^k _n$ is a system of parameters,}
\item{for any $i$, $y^{k+1} _i -y^k _i \equiv 0\; [\MM^{k+1}]$,}
\item{for any $i$, $d(y^k _i)\equiv 0\; [\MM^k]$.}
\end{itemize}
Assume for the moment that such a construction is possible. Then for any $i$,
the sequence $(y^k _i)_{k>0}$ is Cauchy for the $\MM$-adic topology on
$L[[t_1,..,t_n]]$. Since $L[[t_1,..,t_n]]$ is complete, $(y^k _i)_{k>0}$
converges to a formal function $y_i$. By construction, $x_1,y_2,...,y_n$
span the vector space $\MM /\MM^2$, hence it is a system of parameters by
Nakayama Lemma. Moreover by passing to the limit, we find:
$$
\forall i>1, \quad d(y_i)=0
$$
and the result follows. We proceed to this construction by induction on $k>0$.
The case $k=1$ has already been treated above. Assume the construction holds
up to the order $k$. Since $d(y^k _i)\equiv 0\; [\MM^k]$, there exists a
homogeneous polynomial $Q_{i,k}(u_1,...,u_n)$ of degree $k$ such that:
$$
d(y^k _i ) \equiv  Q_{i,k}(x_1,y^1 _2,...,y^1 _n)\; [\MM^{k+1}]
$$
We set $y^{k+1} _i = y^k _i + P_{i,k}(x_1,y^1 _2,...,y^1 _n)$, where every
$P_{i,k}(u_1,...,u_n)$
is defined as:
$$
P_{i,k}(x_1,y^1 _2,...,y^1 _n)=-\int_0 ^{x_1} Q_{i,k}(u,y^1 _2,...,y^1 _n)du
$$
Since every $Q_{i,k}$ is homogeneous of degree $k$, every $P_{i,k}$ is
homogeneous of degree $k+1$. So $y^{k+1} _i -y^k _i \equiv 0\; [\MM^{k+1}]$
for any $i$. In particular, $y^{k+1} _i -y^1 _i \equiv 0\; [\MM^2]$ for any
$i$,
and $x_1,y^{k+1} _2,...,y^{k+1} _n$ span the vector space $\MM/\MM^2$. By
Nakayama Lemma, $x_1,y^{k+1} _2,...,y^{k+1} _n$ is a system of parameters.
Moreover we get by derivation:
$$
d(y^{k+1} _i ) = d(y^{k} _i ) + \frac{\partial P_{i,k}}
{\partial u_1} (x_1,y^1 _2,...,y^1 _n)+
\sum_{i=2} ^n \frac{\partial P_{i,k}}
{\partial u_i} (x_1,y^1 _2,...,y^1 _n)d(y^1 _i)
$$
By construction $\partial P_{i,k}/\partial u_i (x_1,y^1 _2,...,y^1 _n)$ belongs to $\MM^k$ and $d(y^1
_i)$ belongs to $\MM$. By reduction modulo $\MM^{k+1}$ and construction of
$P_{i,k}$, we obtain:
$$
d(y^{k+1} _i ) \equiv d(y^{k} _i ) + \frac{\partial P_{i,k}}
{\partial u_1} (x_1,y^1 _2,..,y^1 _n)\equiv Q_{i,k}(x_1,y^1 _2,.
.,y^1 _n)+\frac{\partial P_{i,k}}
{\partial u_1} (x_1,y^1 _2,..,y^1 _n)\equiv 0\; [\MM^{k+1}]
$$
thus ending the construction to the order $k+1$, and the result follows.
\qed

\begin{prop} \label{formecanon}
Let $d_1,d_2$ be a couple of $L$-derivations on $L[[t_1,..,t_n]]$, where $n>2$.
Assume there
exist two formal functions $x_1,x_2$ such that $d_i(x_i)=1$ for any $i$ and
$d_i(x_j)=0$ for $i\not=j$. Then there exist some formal functions
$y_3,...,y_n$ and $a_3,...,a_n$ such that:
\begin{itemize}
\item{$x_1,x_2,y_3...,y_n$ is a system of parameters,}
\item{$d_1= \frac{\partial}{\partial x_1}$ and $d_2=
\frac{\partial}{\partial x_2} + x_1\sum_{i>2} a_i\frac{\partial}
{\partial y_i} $.}
\end{itemize} 
\end{prop}
\dem Since $d_1(x_1)=1$, lemma \ref{champvecteur} asserts there exist some
formal functions $y' _2,...,y' _n$ such that $x_1,y' _2,...,y' _n$ is a system
of parameters and $d_1(y' _i)=0$ for all $i$. In particular, in this system of
parameters, we have:
$$
d_1 =\frac{\partial}{\partial x_1}
$$ 
A fortiori, this implies that $ker \; d_1 = L[[y' _2,...,y' _n]]$. Since
$d_2(x_1)=0$, there exist some formal functions $b' _2,...,b' _n$ such that:
$$
d_2 = \sum_{i=2} ^n b' _i \frac{\partial}{\partial y' _i}
$$
Set $b_i(y'_2,..,y'_n)=b' _i (0,y'_2,..,y'_n)$ for all $i$, and consider the
$L$-derivation $\partial$:
$$
\partial= \sum_{i=2} ^n b _i \frac{\partial}{\partial y' _i}
$$
By construction, $\partial$ acts on $L[[y' _2,...,y' _n]]$ and
$\partial(f)=d_2(f)(0,y'_2,..,y'_n)$ for any formal function $f$.
Since $d_1(x_2)=0$ and $d_2(x_2)=1$, $x_2$ belongs to
$L[[y' _2,...,y' _n]]$
and $\partial(x_2)=1$. By lemma \ref{champvecteur} applied to
$\partial$ and $L[[y' _2,...,y' _n]]$, there exist
some formal functions $y_3,...,y_n$ in $L[[y' _2,...,y' _n]]$
such that:

\begin{itemize}
\item{$x_2,y_3,...,y_n$ is a system of parameters of $L[[y' _2,...,y' _n]]$, }
\item{for any $i>2$, $\partial (y_i)=0$.}
\end{itemize}
By construction $L[[x_2,y_3,...,y _n]]=L[[y' _2,...,y' _n]]$ and
$x_1,x_2,y_3,...,y_n$ is a system of parameters of $L[[t_1,...,t_n]]$.
Since $ker \; d_1 = L[[y' _2,...,y' _n]]$, this implies:
$$
d_1(x_2)=d_1(y_3)=...=d_1(y_n)=0
$$
In this system of parameters, the derivation $d_1$ can be written as:
$$
d_1=\frac{\partial}{\partial x_1}
$$
Now $d_2(x_1)=0$ and $d_2(x_2)=1$, so that $d_2$ can be written as:
$$
d_2 = \frac{\partial}{\partial x_2}+ \sum_{i>2}
\alpha_i\frac{\partial}{\partial y_i}
$$
where all the $\alpha_i$ are formal functions. Since $d_2(y_i)\equiv
\der(y_i)\equiv 0\;[x_1]$ for any $i$, every function $\alpha_i$ is
divisible by $x_1$. Write $\alpha_i=x_1a_i$, where every $a_i$ belongs to
$L[[t_1,...,t_n]]$. By construction, we find:
$$
d_2 = \frac{\partial}{\partial x_2}+ x_1\sum_{i>2}
a_i\frac{\partial}{\partial y_i}
$$
and the result follows.
\qed

\section{Properties of the derivation $\delta_m$}

In this section, we will analyse the properties of some derivations, and
we will use them to produce the longly-awaited derivation of theorem \ref{principal}. Let $R$ be a commutative
domain with identity. For any positive integer $m$, we define the $R$-derivation
$\delta_m$ on $R[x_1,x_2]$ as:
$$
\delta_m= x_1 ^m \frac{\partial}{\partial x_1} + x_2 ^m \frac{\partial}{\partial
x_2}
$$

\begin{lem} \label{noyau}
For any $m>0$, we have $ker\; \delta_m = R$.
\end{lem}
\dem Note that $\delta_m$ is homogeneous with respect to the variables
$x_1,x_2$. So every element of its kernel can be uniquely written as a sum of
homogeneous elements each belonging to $ker\; \delta_m$. Considering an element
of $ker\; \delta_m$, we may therefore assume that it is homogeneous with respect
to $x_1,x_2$.

Let $P$ be an homogeneous element of $ker\;\delta_m$. If $m=1$, then Euler's
Formula asserts that $\delta_1(P)=deg(P)P$. In this case $deg(P)= 0$ or $P=0$.
So $P$ belongs to $R$ and $ker\; \delta_1 = R$. Assume now that $m>1$, and write
$P$ as:
$$
P(x_1,x_2)=\sum_{k=0} ^n P_k(x_1)x_2 ^k
$$ 
where $P_n(x_1)\not=0$. Then $\delta_m(P)$ can be written as:
$$
\delta_m(P)= \sum_{k=0} ^n kP_k(x_1)x_2 ^{k+m-1} + \sum_{k=0} ^n x_1 ^m P_k '(x_1)x_2 ^k
$$
If $n\not=0$, then the leading term of $\delta_m(P)$ with respect to $x_2$ is
equal to $nP_n(x_1)x_2 ^{n+m-1}$ because $m>1$. In particular
$\delta_m(P)\not=0$ if $n\not=0$. So $n$ must be equal to $0$ and $P=P_0(x_1)$.
But then $\delta_m(P)=x_1 ^m P_0 '(x_1)=0$, which implies that $P_0$ belongs to
$R$, and the result follows.
\qed

\begin{lem} \label{noyau2}
Let $P,Q$ be two elements of $R[x_1,x_2]$, where $Q$ is homogeneous of degree
$k$ with respect to $x_1,x_2$. Assume that they satisfy the following relation:
$$
\delta_m(P)+ x_1x_2 ^m Q=0
$$
If $m\geq k+4$, then $P$ belongs to $R$ and $Q=0$.
\end{lem}
\dem Up to adding a constant to $P$, we may assume that $P(0,0)=0$. We are going
to prove that $P=Q=0$. By the
previous lemma, $\delta_m$ is injective on the polynomials with no constant
terms. Since $Q$ and $\delta_m$ are homogeneous, $P$ needs to be homogeneous.
We write $P$ as the sum:
$$
P(x_1,x_2)= \sum_{i+j=r} a_{i,j} x_1 ^i x_2 ^j 
$$
By derivation, we find:
$$
\delta_m(P)=\sum_{i+j=r} ia_{i,j} x_1 ^{i+m-1} x_2 ^j + \sum_{i+j=r} ja_{i,j} x_1 ^{i}
x_2 ^{j+m-1}= -x_1x_2 ^m Q(x_1,x_2)
$$
Since $Q$ is homogeneous of degree $k$, we have $r-2=k$ and $r\leq m-2$, since
$m\geq k+4$. We have $m-1>0$ and the following congruence holds
modulo $x_2 ^{m-1}$:
$$
\delta_m(P)\equiv \sum_{i+j=r} ia_{i,j} x_1 ^{i+m-1} x_2 ^j \equiv 0\; [x_2
^{m-1}]
$$ 
In this sum, all the indices $j$ in this sum satisfy $\leq r\leq m-2$, so that
$ia_{i,j}=0$ for every couple $(i,j)$. Therefore $a_{i,j}=0$ or $i=0$, and $P$
reduces to a polynomial of the form $ax_2 ^r$, where $a$ belongs to $R$.
But then we find:
$$
\delta_m(P)= a r x_2 ^{m+r-1} = -x_1 x_2 ^{m} Q(x_1,x_2)
$$
which is impossible unless $P=Q=0$.
\qed 

\section{Proof of the main theorem}

In this section, we work under the assumptions of section \ref{redd}. Let
$A$ be an integral $k$-algebra of finite type over a field of characteristic zero.
Let $d_1,d_2$ be a couple of $k$-derivations on $A$, and $x_1,x_2$ be two
elements of $A$ such that $d_i(x_i)=1$ for any $i$ and $d_i(x_j)=0$ for
$i\not=j$. We assume that $k=ker\; d_1 \cap ker\; d_2$. According to
section \ref{redd}, $A$ is embedded in $L[[t_1,...,t_n]]$, where $L$
is a finite extension of
$k$, and the $d_i$ extend to $L$-derivations on $L[[t_1,...,t_n]]$.
By proposition \ref{formecanon}, there exist some formal elements
$y_3,..,y_n,a_3,..,a_n$ of $L[[t_1,...,t_n]]$ such that:

\begin{itemize}
\item{$x_1,x_2,y_3...,y_n$ is a system of parameters,}
\item{$d_1= \frac{\partial}{\partial x_1}$ and $d_2=
\frac{\partial}{\partial x_2} + x_1\sum_{i>2} a_i\frac{\partial}
{\partial y_i} $.}
\end{itemize} 
For any positive integer $m$, we consider the following $k$-derivation
$\Delta_m$:
$$
\Delta_m = x_1 ^m d_1 + x_2 ^m d_2
$$
In order to prove theorem \ref{principal}, we are going to establish
the following result.

\begin{prop} \label{demons}
Under the previous assumptions, there exists a positive integer $m_0$
such that, for any $m\geq m_0$, we have $ker\; \Delta_m =ker\; d_1 \cap
ker\; d_2=k$ in $A$.
\end{prop}
Before giving the proof, we make explicit the number $m_0$. Let $R$ be the
ring $L[[y_3,...,y_n]]$. We provide $R[[x_1,x_2]]=L[[t_1,..,t_n]]$
with the homogeneous degree $deg$ on the variables $x_1,x_2$. More precisely
$deg(x_1)=deg(x_2)=1$ and $deg(y_i)=0$ for any $i\geq 3$. Let $a^k _i$ be the
homogeneous part of $a_i$ of degree $k$. We denote by $\partial_k$ the
$L$-derivation on $L[[t_1,..,t_n]]$:
$$
\der_k = \sum_{i>2} a^k _i\frac{\partial} {\partial y_i}
$$
Every $\der_k$ is an $L$-linear operator of degree $k$ on $R[[x_1,x_2]]$, and
we have by construction:
$$
d_2=\frac{\partial}{\partial x_2} + x_1\sum_{k\geq 0} \der_k
$$
Let ${\cal{D}}=Der_L(L[[t_1,..,t_n]],L[[t_1,..,t_n]])$ be the space of
$L$-derivations
on $L[[t_1,..,t_n]]$, and let $M$ be the sub-$R[[x_1,x_2]]$-module of
${\cal{D}}$ spanned by the $\der_k$. Since ${\cal{D}}$ is noetherian, there
exists an integer $m_0 '$ such that:
$$
M= R[[x_1,x_2]]\left\{\der_0 ,...,\der_{m_0 '}\right\}
$$
We denote by $m_0$ the integer $m_0 = m_0 '+4$. \\ \\
{\it Proof of proposition \ref{demons}}: Let $f$ be a non-zero element of
$R[[x_1,x_2]]$ such that $\Delta_m(f)=0$, and assume that $m\geq m_0$.
Let us prove by contradiction that $d_1(f)=d_2(f)=0$.
We decompose $f$ into its sum of homogeneous parts with respect to $deg$:
$$
f=f_0 +f_1 +...+f_k+...
$$
{\underline {Assertion 1}:} $f_0\not=0$. \\ \\
Let $i$ be the smallest index such that $f_i\not=0$. Since $\Delta_m(f)=0$,
we have $\delta_m(f) + x_1x_2 ^m\sum_{k\geq 0} \der_k(f)=0$.
By considering only the terms of degree $i+m-1$ in this equation, we find that
$\delta_m(f_i)=0$. By lemma \ref{noyau}, $f_i$ belongs to $R$ and $i=0$. \\ \\
{\underline {Assertion 2}:} $f\not=f_0$. \\ \\
Assume on the contrary that $f=f_0$. Since $d_1=\der/\der x_1$, we have
$d_1(f)=0$. As $\Delta_m= x_1 ^m d_1 + x_2 ^m d_2$ and $\Delta_m(f)=0$,
we get $d_1(f)=d_2(f)=0$, hence a contradiction. \\ \\
{\underline {Assertion 3}:} Let $k$ be the smallest positive index such that
$f_k\not=0$. Then $2\leq k\leq m_0 ' +2$. \\ \\
If we consider the terms of degree $m+i-1$ in the equation $\Delta_m(f)=0$, we
obtain:
$$
\delta_m(f_i) + x_1x_2 ^m \left\{ \der_0(f_{i-2})+ \der_1(f_{i-3})+
...+
\der_{i-2}(f_0)\right\}=0
$$
for any $i>1$, and $\delta_m(f_1)=0$. By lemma \ref{noyau}, we have $f_1=0$
and $k\geq 2$. Assume on the contrary that $k>m_0 ' +2$.
Then for any $0<i<k$, we have $f_i=0$. This implies that $\der_{i}(f_0)=0$
for any $i<k-2$, and we have in particular:
$$
\der_0(f_0)=\der_1(f_0)=...=\der_{m_0 '}(f_0)=0
$$
By assumption, $\der_0 ,...,\der_{m_0 '}$ span the module $M$ generated by the
$\der_l$. For any nonnegative integer $l$, there exist some formal elements
$b_{l,\beta}$ such that:
$$
\der_l=\sum_{\beta=0} ^{m_0 '} b_{l,\beta} \der_{\beta}
$$
This implies that $\der_l(f_0)=0$ for any $l$. In particular $f_0$ is
annihilated by $d_2$. Since $d_1(f_0)=0$, $f_0$ belongs to $ker\; \Delta_m$.
So $(f-f_0)$ lies in $ker\; \Delta_m$ and has no homogeneous part of degree
$0$. By our first assertion, $(f-f_0)$ needs to be equal to $0$. But that
contradicts our second assertion. \\ \\
{\underline {Assertion 4}:} The final contradiction. \\ \\
With $k$ as in Assertion 3, if we consider the terms of degree $m+k-1$ in the equation $\Delta_m(f)=0$, we
obtain:
$$
\delta_m(f_k) + x_1 x_2 ^m \der_{k-2} (f_0)=0
$$
because $f_1=...=f_{k-1}=0$. By the previous assertion, $k\leq m_0 ' +2$ and
$\der_{k-2} (f_0)$ is an homogeneous polynomial of degree $\leq m_0 '$. In
particular, we find:
$$
m\geq m_0 \geq m_0 ' +4 \geq deg(\der_{k-2} (f_0))+4
$$
By applying lemma \ref{noyau2} to the relation $\delta_m(f_k) + x_1 x_2 ^m
\der_{k-2} (f_0)=0$, we get that $f_k=0$, which is
impossible
by the third assertion.
\qed
Theorem \ref{principal} is a direct consequence of this proposition. Indeed
if $m\geq m_0$ and $f$ is an element of $A$ such that $\Delta_m(f)=0$, then
$d_1(f)=d_2(f)=0$ as elements of $L[[t_1,..,t_n]]$. Since $A$ is a subring of
$L[[t_1,..,t_n]]$, $f$ belongs to $k$ and the result follows.

\section{A few consequences and an example}

In this section, we are going to derive some consequences of theorem
\ref{principal}. In particular we will give another proof of Nowicki's Theorem.

\subsection{Proof of theorem \ref{cons2}}

Let $k$ be an algebraically closed field of characteristic zero. Let $A$ be a
$k$-algebra of finite type with no zero divisors. Let $\FF$ be a family of
$k$-derivations on $A$, and denote by $M_{\FF}$ the $A$-module spanned by $\FF$.
Recall that $M_{\FF,min}$ is the set of $\FF$-minimal derivations, i.e. the
subset of $M_{\FF}$ formed by the $k$-derivations $d$ such that:
$$
ker\; \FF = ker \; d
$$
Let us prove that $M_{\FF,min}$ is a non-empty residual subset of $M_{\FF}$.
The non-emptyness is guaranteed by theorem \ref{principal}. For residuality,
let $F$ be any finite $k$-vector subspace of $M_{\FF}$. Since $A$ is a finite
$k$-algebra, its dimension is at most countable. So the space $A/ker\; \FF$
admits a filtration $\{F_n\}_{n\in \mathbb{N}}$ where every $F_n$ is finite
dimensionnal, i.e.
$$
\frac{A}{ker\; \FF} = \cup_{n\in \mathbb{N}} \; F_n \quad \mbox{and} \quad
F_0 \subseteq F_1 \subseteq ...\subseteq F_n \subseteq ...
$$
Consider the set $\Sigma_n$ defined as:
$$
\Sigma_n = \left\{(d,f) \in F \times \mathbb{P}(F_n), \; d(f)=0\right \}
$$
Note that $\Sigma$ is an algebraic subset of $F \times \mathbb{P}(F_n)$
by construction. If $\Pi:F \times \mathbb{P}(F_n)\rightarrow F$ denotes
the standard projection on $F$, we find:
$$
\Pi(\Sigma_n) = \left\{d \in F, \; \exists d \in \mathbb{P}(F_n), \; d(f)=0\right \}
$$
In particular, a $k$-derivation $d$ in $F$ is $\FF$-minimal if and only if it
does not belong to the union $\cup_{n\in \mathbb{N}} \; \Pi(\Sigma_n)$. In other
words we have:
$$
M_{\FF, min} \cap F = F - \cup_{n\in \mathbb{N}} \; \Pi(\Sigma_n)=
\cap_{n\in \mathbb{N}} \; (F-\Pi(\Sigma_n))
$$
But every $\Sigma_n$ is closed in $F \times \mathbb{P}(F_n)$ for the Zariski topology.
Since $\mathbb{P}(F_n)$ is projective, it is a complete variety and the projection $\Pi$
is closed. Therefore every set $\Pi(\Sigma_n)$ is closed in $F$, and $M_{\FF,min}$
intersects $F$ in a countable intersection of Zariski open sets. By definition,
$M_{\FF.min}$ is residual in $M_{\FF}$.

\subsection{Proof of Nowicki's result}

In this subsection, we will use theorem \ref{principal} in order to establish
theorem \ref{Nowicki2}. The proof is based on the following lemma.

\begin{lem} \label{extension}
Let $A$ be a finite $k$-algebra with no-zero divisors. Let $B$ be a sub-algebra
of $A$ and let $d_B$ be a $k$-derivation on $B$. Then there exists a non-zero
element $a$ of $A$ and a $k$-derivation $d_A$ on $A$ such that $d_A =a d_B$
on $B$.
\end{lem}
\dem Let $x_1,...,x_n$ be a set of generators of $A$ as a $k$-algebra. If $K=Fr(B)$
and $L=Fr(A)$, then we have $L=K(x_1,...,x_n)$. Consider the following chain of fields:
$$
K=K_0 \subseteq K_1 \subseteq ...\subseteq K_n \quad \mbox{where} \quad
K_i=K(x_1,...,x_i)
$$
By construction we have $K_{i+1}=K_i(x_{i+1})$. Assume there exists a $k$-derivation
$\der_i$ on $K_i$. Then there exists a $k$-derivation $\der_{i+1}$ on $K_{i+1}$
that extends the derivation $\der_i$ (if $x_i$ is transcendental over $K_i$, see
\cite{Ka}, p.10 and if $x_i$ is algebraic over $K_i$, see \cite{Ma}, p. 14). By a finite
induction, we can extend the $k$-derivation $d_B=\der_0$ of $K=K_0$ into a $k$-derivation
$d_L=\der_n$ on $L=K_n$. For any index $i$, there exists a non-zero element $a_i$ of $A$
such that $a_id_L (x_i)$ belongs to $A$. If $a=a_1 ...a_n$, then $a d_L(x)$ belongs to
$A$ for any $x$ in $A$. The derivation $d_A=a d_L$ maps $A$ into $A$ and thus defines a
$k$-derivation on $A$. By construction we have $d_A = a d_B$ on $B$.
\qed
{\it Proof of theorem \ref{Nowicki2}}: Let $A$ be a finite $k$-algebra with no zero
divisors over a field $k$ of characteristic zero. Let $d$ be a $k$-derivation on
$A$ and set $B=ker \; d$. By Leibniz rule, $B$ needs to be a $k$-algebra. Let
$a$ be any element of $A$ that is algebraic over $B$. Let $P(t)=b_0 + ...+ b_n t^n$
be a polynomial of $B[t]$ of minimal degree such that $P(a)=0$. By derivation,
we get:
$$
d(P(a))= \left(\sum_{k=1} ^n ka^{k-1}\right)d(a)= P'(a)d(A)=0
$$
By minimality of the degree, we find $P'(a)\not=0$ and $d(A)=0$. So $a$ belongs to $B$
and $B$ is algebraically closed in $A$.

Conversely let $B$ be a sub-algebra of $A$ and assume that $B$ is algebraically closed
in $A$. For any $x$ outside $B$, denote by $B_x$ the algebra $B[x]$ and consider the
$k$-derivation $d_x=d/dx$ on $B_x$. Since $x$ lies outside $B$, $x$ is transcendental
over $B$ and $d_x$ is well-defined. By lemma \ref{extension}, there exist a
a non-zero element $a_x$ of $A$ and a $k$-derivation $\der_x$ on $A$ such that
$\der_x =a_x d_x$. In particular $\der_x(x)=a_xd_x(x)=a_x\not=0$. Consider the family
$\FF=\{\der_x\}_{x\in A-B}$ of $k$-derivations on $A$. By construction, we have:
$$
ker\; \FF=B
$$
By theorem \ref{principal}, there exists an $A$-linear combination $d$ of the $\der_x$
such that $ker\; d = B$. In particular $B$ is the kernel of a $k$-derivation and the
result follows.
\qed

\subsection{An example}

In this subsection, we will illustrate the notion of residuality in an example.
Consider the algebra $A=\CC[x,y]$ provided with the following
couple of $\CC$-derivations:
$$
\FF=\left\{d_1,d_2\right\}=\left\{ x\frac{\der }{\der x}, y\frac{\der}{\der y} \right\}
$$
Let $F$ be the vector space spanned by $d_1,d_2$. Then it is easy to check that
an element $\lambda_1 d_1 + \lambda_2 d_2$ of $F$ admits a first integral if and
only if either $\lambda_2=0$ or $\lambda_1/\lambda_2$ is a nonpositive
rational number. In this latter case, if $\lambda_1/\lambda_2=-p/q$ where
$p\geq 0,q>0$ are coprime integers, then every first integral of 
$\lambda_1 d_1 + \lambda_2 d_2$ is a polynomial function of the
expression:
$$
f_{p,q}(x,y)= x^q y^p
$$
Let ${\cal{D}}$ be the union of the coordinate axes of $\CC^2$ and of the
lines of rational negative slopes. Then ${\cal{D}}$ is a countable union
of complex lines in $\CC^2$. If we identify $d_1,d_2$ with the canonical
basis of $\CC^2$, then we have:
$$
F \cap M_{\FF,min} = F-{\cal{D}}
$$
and this latter set turns out to be residual, since it is a countable intersection
of Zariski open sets of $F\simeq \CC^2$. Residuality occurs as we cannot find
a finite dimensionnal space in $\CC[x,y]$ containing all
the first integrals of minimal degree of elements of $F$.
If it were possible, then as in theorem \ref{cons2}, we could prove
that $F \cap M_{\FF,min}$ would be a finite
intersection of Zariski open sets, hence an open set and this is not the case.


\begin{thebibliography}{Da-F}

\bibitem[Da-F]{Da-F} D.Daigle, G.Freudenburg {\it A counterexample to Hilbert's fourteenth problem in dimension 5}, J. Algebra
221 (1999), $n^0 2$ pp. 528-535.

\bibitem[De-F]{De-F} J.Deveney, D.Finston {\it $G_a$-actions on $\CC^3$ and $\CC^7$}, Communications in Algebra, 22(15), pp.6295-6302
(1994).

\bibitem[Ei]{Ei} D.Eisenbud {\it Commutative Algebra with a view toward Algebraic
Geometry}, Springer Verlag New York (1995).

\bibitem[Hum]{Hum} J.Humphreys {\it Linear algebraic groups}, Graduate Texts in Math. $n^0 21$, Springer Verlag New York Heidelberg,
1975.

\bibitem[Jou]{Jou} J-P.Jouanolou {\it Equations de Pfaff alg\'ebriques}, Lect. Notes in Math. 708, Springer Verlag Berlin (1979).

\bibitem[Ka]{Ka} I.Kaplansky {\it An introduction to Differential Algebra}, publications de l'institut
de math\'ematiques de Nancago, Hermann (1957).

\bibitem[Ma]{Ma} A.R.Magid {\it Lectures on Differential Galois Theory}, University Lecture Series vol.7,
American Mathematical Society (1994).

\bibitem[Na]{Na} M. Nagata {\it On the fourteenth problem of Hilbert}, Proc. Intern. Congress Math. (1958), pp. 459-462, Cambridge
University Press, New York 1966.

\bibitem[Na2]{Na2} M.Nagata, A.Nowicki {\it Rings of constants for $k$-derivations in $k[x_1,...,x_n]$}, J. Math. Kyoto Univ. 28
(1998) $n^0 1$, pp. 111-118.

\bibitem[No]{No} A.Nowicki {\it Rings and fields of constants for derivations
in characteristic zero}, J. Pure Appl. Algebra 96 (1994), $n^0 1$, pp. 47-55.

\bibitem[No2]{No2} A.Nowicki {\it On the nonexistence of rational first integrals for systems of linear differential equations},
Linear Algebra Appl. 235 (1996), pp. 107-120.


\end{thebibliography}
\end{document}